# Symbolic estimation of distances between eigenvalues of Hermitian operator

Ilia Lomidze, Natela Chachava (Georgian Technical university, Tbilisi, Georgia)
E-mail: lomiltsu@gmail.com,  chachava.natela@yahoo.com

Abstract
   We find out a method for symbolic estimation of a minimal (maximal) distance between eigenvalues of a Hermitian matrix (or roots of a polynomial with real (maybe degenerated) roots), using Hankel matrices formalism. The range of location of eigenvalues is symbolically estimated too. All estimations can be done with any precision.

Introduction

Localization of roots of polynomial was considered with a number of authors:  Jacoby has defined the number of roots and the number of real roots applying  Hankel matrix method and received the appropriate formula; Joachimsthal  defined the number of (real) roots lying in certain interval using the Sturm consequence; Gershgorin and Cauchy considered the matrixes' (maybe complex) eigenvalues' localization, this method may be clarified by some numerical methods. We restricted ourselves with polynomials with *real* roots (our aim is to describe states of some quantum system by (Hermitian) density matrix), but have considered *degenerated* roots too. The main method we have used was proposed in [1]. Similarly to the Jacobi theorem, it allows us to define the number of different roots of Hermitian matrix. We have used formulae for main minors [1] to express degree of degeneracy of the roots (or eigenvalues). Using the method proposed we derived the formulae for the multiplicity of the certain root and for coefficients of minimal polynomial. In the second part algorithm constructed is applied to a polynomial of a special kind which allows us to build a recursive sequences for minimal and maximal distances between roots of given polynomial. The limits of these sequences are exact values of minimal and maximal distances. The members of sequences are *rational functions* of the coefficients of the polynomial. This method also allows us to find the segment which includes all roots expressing its ends as rational functions of the coefficients of the given polynomial with any precision. In the last part, as illustration of the method developed, it is applied to Wilkinson's type polynomial.

Let's consider polynomial $P_n(x)$ with real roots $p_1 < p_2 < ... < p_m$ having multiplicities $r_1, r_2, ..., r_m$ ($\sum_{i=1}^{m} r_i = n$). Some notations must be introduced (here we are following to the article [1]):

$$t_k = \sum_{l=1}^{m} r_l p_l^k, \quad k = 0, 1, 2, ..., \qquad (1)$$

$$H_k = \left[ h_{ij} \right]_1^k = \left[ t_{i+j-2} \right]_1^k, \qquad (2)$$

$H_k$ is *k*-order Hankel matrix [2]. Taking into account that

$$t_{i+j-2} = \sum_{l=1}^{m} p_l^{i-1} r_l p_l^{j-1} = \left[ p_l^{i-1} \right]^T \left[ r_l p_l^{j-1} \right],$$

($\left[ p_l^{i-1} \right]$ denotes one-column matrix, and $\left[ p_l^{i-1} \right]^T$ - correspondent one-row matrix) and using the Binet-Cauchy Theorem [3] it is easy to show (see [1]) that

$$D_k = \det H_k = \begin{cases} \sum_{1 \le i_1 < \cdots < i_k \le m} \left\{ r_{i_1} \cdots r_{i_k} \prod_{1 \le j < l \le k} (p_{i_j} - p_{i_l})^2 \right\} > 0, & 1 \le k \le m, \\ 0, & k > m. \end{cases} \qquad (3)$$

In particular, one has

$$D_2 = \det \begin{bmatrix} t_0 & t_1 \\ t_1 & t_2 \end{bmatrix} = \det \left( \begin{bmatrix} 1 & \cdots & 1 \\ p_1 & \cdots & p_m \end{bmatrix} \begin{bmatrix} r_1 & \cdots & r_m \\ r_1 p_1 & \cdots & r_m p_m \end{bmatrix}^T \right) = \sum_{1 \le i < j \le m} r_i r_j (p_j - p_i)^2, \qquad (3')$$

$$D_m = \det H_m = \left( \prod_{i=1}^m r_i \right) \prod_{1 \le i < j \le m} (p_j - p_i)^2. \qquad (3'')$$

Note, that all sums (1) and therefore all determinants (3) can be rationally expressed via coefficients of the polynomial, using known Newton's formulae.

In the article [1] it was shown that any $m$ equalities from (1) can be solved for $r_i$, $i = \overline{1, m}$:

$$r_i^{-1} = D_m^{-1} \det \begin{bmatrix} H_m & (\mathbf{p}_i)^T \\ -\mathbf{p}_i & 0 \end{bmatrix}, \qquad (4)$$

where

$$\mathbf{p}_i \equiv \left[ p_i^{k-1} \right]_{k=\overline{1,m}} = \left[ 1, p_i, \ldots, p_i^{m-1} \right]$$

is a row, whereas $\mathbf{p}_i^T$ is a correspondent column.

Calculating the determinant in (4), we obtain

$$r_i^{-1} = \langle \mathbf{p}_i H_m^{-1} \mathbf{p}_i^T \rangle, \qquad (5)$$

where $H_m^{-1}$ denotes the inverse matrix and brackets stand for a scalar product of a row by a column.

**Statement 1.** The formula (5) can be generalized as:

$$r_i^{-1} \delta_{ij} = \langle \mathbf{p}_i H_m^{-1} \mathbf{p}_j^T \rangle, \qquad (6)$$

where $\delta_{ij}$ denotes the Kronecker delta.

**Proof.** It is enough to present the inverse matrix $H_m^{-1}$ as

$$H_m^{-1} = \left( \left[ p_l^{i-1} \right] \left[ r_l p_l^{j-1} \right]^T \right)^{-1} = \left( \left[ r_l p_l^{j-1} \right]^{-1} \right)^T \left[ p_l^{i-1} \right]^{-1}$$

and then calculate the scalar product (6).

**Corollary 1.** The next formula is valid

$$r_j^{-1} - r_i^{-1} = \langle (\mathbf{p}_j + \mathbf{p}_i) H_m^{-1} (\mathbf{p}_j - \mathbf{p}_i)^T \rangle. \qquad (7)$$

**Proof.** Obvious.

Expanding the difference of degrees in the last column in the formula (7) we obtain

$$p_j^{k-1} - p_i^{k-1} = (p_j - p_i) \sum_{l=0}^{k-2} p_i^{k-2-l} p_j^l, \quad k = \overline{2, m}. \qquad (8)$$

Let's construct the expression

$$\frac{r_j^{-1} - r_i^{-1}}{p_j - p_i} \equiv \Delta_{ij} = \Delta_{ji} = \sum_{n,k=1}^m (p_j^{n-1} + p_i^{n-1}) [H_m^{-1}]_{nk} \sum_{l=0}^{k-2} p_i^{k-2-l} p_j^l \qquad (9)$$

(In the right hand side of the last formula the summand with $k=1$ is equal to zero). Sum of all $\Delta_{ij}$ ($i, j$, $1 \le i < j \le m$), obviously, gives us the symmetric function of $p_k$ in the right hand side of the (9).

The symmetric function can be expressed via the elementary symmetric functions of the variable $p_k$

$$\sigma_1 = \sum_{1 \leq j \leq m} p_j, \quad \sigma_2 = \sum_{1 \leq j < k \leq m} p_j p_k, \ldots, \quad \sigma_m = p_1 p_2 \cdots p_m, \tag{10}$$

(or, alternatively, via sums $s_k = \sum_{i=1}^{m} p_i^{k-1}$, $k = \overline{1, m}$, using the Newton's formulae).

For forthcoming consideration we need the next theorem:

**Theorem 1.** All elementary symmetric functions (10) can be rationally expressed via sums (1) and therefore via the coefficients of the given polynomial $P_n(x)$.

**Proof.** Let us consider the polynomial

$$P_m(x) = D_m^{-1} \det \begin{bmatrix} t_0 & t_1 & \cdots & t_{m-1} & t_m \\ t_1 & t_2 & \cdots & t_m & t_{m+1} \\ \cdots & \cdots & \cdots & \cdots & \cdots \\ t_{m-1} & t_m & \cdots & t_{2m-2} & t_{2m-1} \\ 1 & x & \cdots & x^{m-1} & x^m \end{bmatrix} = D_m^{-1} \det \begin{bmatrix} H_m & (\mathbf{t})^T \\ \mathbf{x} & x^m \end{bmatrix}, \tag{11}$$

having the first coefficient 1. Let us show that other coefficients of this polynomial are from the set (10).

Let us use the obvious relation

$$\begin{bmatrix} H_m & (\mathbf{t})^T \\ \mathbf{x} & x^m \end{bmatrix} = \begin{bmatrix} [p_l^{i-1}]_1^m & \mathbf{0}^T \\ \mathbf{0} & 1 \end{bmatrix} \left( \begin{bmatrix} r_l p_l^{j-1} \end{bmatrix}_1^m \right)^T \begin{bmatrix} r_l p_l^m \end{bmatrix}^T .$$

Then, taking into account the formulae (3), (11) and the Viet theorem one easily finds

$$P_m(x) = D_m^{-1} \det [p_i^{k-1}]_1^m \left( \prod_{i=1}^m r_i \right) \det \begin{bmatrix} \left( [p_l^{j-1}]_1^m \right)^T & [p_l^m]^T \\ \mathbf{x} & x^m \end{bmatrix} =$$

$$= D_m^{-1} \left( \prod_{i=1}^m r_i \right) \left( \det [p_i^{k-1}] \right)^2 \prod_{i=1}^m (x - p_i) = \prod_{i=1}^m (x - p_i) =$$

$$= x^m + \sum_{k=1}^m x^{m-k} (-1)^k \sigma_k.$$

Hence, one gets

$$\sigma_k = \sum_{1 \leq i < j \leq m} p_{i_1} \cdots p_{i_k} = D_m^{-1} H \begin{pmatrix} m+1 \\ m-k+1 \end{pmatrix}, \tag{12}$$

where $H \begin{pmatrix} m+1 \\ m-k+1 \end{pmatrix}$ denotes the minor which remains after deleting the $(m+1)$- th row and the $(m-k+1)$- th column, $k = \overline{1, m}$, of the determinant in the formula (11). These minors, as well as the factor $D_m^{-1}$, are polynomial functions of the sums (1) and therefore the statement of the Theorem 1 follows. ∎

**Corollary 2.** The polynomial (11) is a minimal polynomial of Hermitian operator (matrix) having the polynomial $P_n(x)$ as a characteristic one.

**Corollary 3.** The function

$$Z \equiv \sum_{i<j} \frac{r_j^{-1} - r_i^{-1}}{p_j - p_i} = \sum_{i<j} \left\{ \sum_{n,k=1}^{m} (p_j^{n-1} + p_i^{n-1})[H_m^{-1}]_{nk} \sum_{l=0}^{k-2} p_i^{k-2-l} p_j^l \right\} \qquad (13)$$

is a rational function of the coefficients of the polynomial $P_n(x)$.

Let us apply the formula (13) to the polynomial

$$Q_{3m}(x) = P_m(x) P_m^2(x - \varepsilon) = \prod_{i=1}^{m}(x - p_i) \prod_{j=1}^{m}(x - (p_j + \varepsilon))^2, \qquad (14)$$

assuming that the $\varepsilon$ is chosen according to the condition

$$p_i - (p_j + \varepsilon) \neq 0, \qquad i, j = \overline{1, m}, \ m \geq 3. \qquad (15)$$

After obvious transformations one gets

$$Z(\varepsilon) = \left(1 - \frac{1}{2}\right) \sum_{i,j=1}^{m} \frac{1}{p_i - (p_j + \varepsilon)} = \frac{1}{2} \left( \sum_{i<j} \frac{1}{p_i - (p_j + \varepsilon)} - \frac{m}{\varepsilon} + \sum_{i>j} \frac{1}{p_i - (p_j + \varepsilon)} \right) =$$

$$= \frac{1}{2} \left( -\frac{m}{\varepsilon} + \sum_{1 \leq i < j \leq m} \frac{2\varepsilon}{(p_i - p_j)^2 - \varepsilon^2} \right). \qquad (16)$$

$$(m \geq 3)$$

Let us assume that

$$\varepsilon = \mu_0, \quad 0 < \mu_0 < \min |p_i - p_j| \equiv \mu.$$

Then one obtains

$$\frac{1}{\mu_0} \left( Z(\mu_0) + \frac{m}{2\mu_0} \right) = \sum_{1 \leq i < j \leq m} \frac{1}{(p_i - p_j)^2 - \mu_0^2} > \frac{1}{\mu^2 - \mu_0^2} > 0. \qquad (17)$$

As result, the next estimation follows to the (17)

$$\mu^2 > \left( \sum_{1 \leq i < j \leq m} \frac{1}{(p_i - p_j)^2 - \mu_0^2} \right)^{-1} + \mu_0^2 \equiv \mu_1^2. \qquad (18)$$

Obviously,

$$\mu_0^2 < \mu_1^2 < \mu^2.$$

Hence, one can use $\mu_1$ instead $\mu_0$ in the formula (16) to calculate $\mu_2 > 0$:

$$\mu_2^2 \equiv \left( \sum_{1 \leq i < j \leq m} \frac{1}{(p_i - p_j)^2 - \mu_1^2} \right)^{-1} + \mu_1^2 < \mu^2.$$

For similar reasons one gets

$$\mu_0^2 < \mu_1^2 < \mu_2^2 < \mu^2.$$

Continuing this process, after $k$ steps one obtains the increasing and over bounded sequence

$$\mu_0^2 < \mu_1^2 < \mu_2^2 < ... < \mu_k^2 < \mu^2, \qquad (19)$$

where

$$\mu_k^2 \equiv \left( \sum_{1 \leq i < j \leq m} \frac{1}{(p_i - p_j)^2 - \mu_{k-1}^2} \right)^{-1} + \mu_{k-1}^2, \quad k = 1, 2, ... \ . \qquad (20)$$

Therefore, when $k \to \infty$ the sequence (19) must converge to some limit $\lim_{k \to \infty} \mu_k = \mu_\infty$, $0 < \mu_\infty \leq \mu$. This limit can be easily calculated using the recurrent relation (20): assuming $k \to \infty$ in the both sides of it, one obtains

$$\mu_\infty^2 \equiv \left( \sum_{1 \leq i < j \leq m} \frac{1}{(p_i - p_j)^2 - \mu_\infty^2} \right)^{-1} + \mu_\infty^2.$$

Hence,

$$\left( \sum_{1 \leq i < j \leq m} \frac{1}{(p_i - p_j)^2 - \mu_\infty^2} \right)^{-1} = 0, \qquad \sum_{1 \leq i < j \leq m} \frac{1}{(p_i - p_j)^2 - \mu_\infty^2} \to \infty.$$

So, in the last sum at least one summand is infinite – at least one denominator in it is equal to 0. Taking into account that $0 < \mu_\infty \leq \mu \equiv \min|p_i - p_j|$, one has to conclude

$$\mu_\infty = \mu = \min|p_i - p_j|. \tag{21}$$

Similarly, assuming in the formula (16)

$$\varepsilon = M_0 > \max|p_i - p_j| = p_m - p_1 \equiv M > 0,$$

one gets from the formula (17):

$$\frac{1}{M_0}\left( Z(M_0) + \frac{m}{2M_0} \right) = \sum_{1 \leq i < j \leq m} \frac{1}{(p_i - p_j)^2 - M_0^2} < \frac{1}{M^2 - M_0^2} < 0. \tag{22}$$

Taking into account that all summands in the formula (22) are negative, after obvious transformations, one obtains:

$$M^2 < \left( \sum_{1 \leq i < j \leq m} \frac{1}{(p_i - p_j)^2 - M_0^2} \right)^{-1} + M_0^2 \equiv M_1^2,$$

$$M^2 < M_1^2 < M_0^2.$$

As far as $M_1^2 > M^2 = \max(p_i - p_j)^2$, one can use $M_1^2$ instead $M_0^2$ in the formula (22) to calculate $M_2^2$:

$$M^2 < \left( \sum_{1 \leq i < j \leq m} \frac{1}{(p_i - p_j)^2 - M_0^2} \right)^{-1} + M_1^2 \equiv M_2^2,$$

$$M^2 < M_2^2 < M_1^2 < M_0^2.$$

Continuing this process, after *k* steps one obtains the decreasing and lower bounded sequence:

$$M^2 < M_k^2 < \ldots < M_2^2 < M_1^2 < M_0^2, \tag{23}$$

where

$$M^2 < \left( \sum_{1 \leq i < j \leq m} \frac{1}{(p_i - p_j)^2 - M_{k-1}^2} \right)^{-1} + M_{k-1}^2 = M_k^2. \tag{24}$$

Therefore, the sequence (23) must converge to some limit $\lim_{k \to \infty} M_k = M_\infty$, $0 < M_\infty \leq M$, when $k \to \infty$. This limit can be easily calculated using the recurrent relation (24): assuming $k \to \infty$ in the both sides of it, one obtains

$$\left( \sum_{1 \leq i < j \leq m} \frac{1}{(p_i - p_j)^2 - M_\infty^2} \right)^{-1} + M_\infty^2 = M_\infty^2. \tag{25}$$

This means that

$$\left( \sum_{1 \leq i < j \leq m} \frac{1}{(p_i - p_j)^2 - M_\infty^2} \right)^{-1} = 0, \qquad \sum_{1 \leq i < j \leq m} \frac{1}{(p_i - p_j)^2 - M_\infty^2} \to \infty. \tag{26}$$

So, in the last sum at least one summand is infinite – at least one denominator in it is equal to 0. Taking into account that $0 < \max |p_i - p_j| = M \leq M_\infty$ one has to conclude

$$M_\infty = M = \max |p_i - p_j|. \qquad (27)$$

In order to perform this scheme really one has to prove the next

**Lemma**: The parameters $\mu_k$, $0 < \mu_k < \mu$ and $M_k$, $0 < M_k < M$, $k = 0, 1, 2, \ldots$, can be expressed through known coefficients of the polynomial $P_n(x)$.

**Proof of the Lemma.** One can observe from the formula (17) that

$$\frac{1}{\mu_0}\left(Z(\mu_0) + \frac{m}{2\mu_0}\right)\bigg|_{\mu_0 \to 0+} = \sum_{1 \leq i < j \leq m}(p_i - p_j)^{-2} > \frac{1}{\mu^2} > 0. \qquad (28)$$

Hence,

$$\mu^2 > \left[\sum_{1 \leq i < j \leq m}(p_i - p_j)^{-2}\right]^{-1}. \qquad (29)$$

Therefore, one can start the process using as initial $\mu_0 > 0$ the value defined by the formula (29):

$$\mu_0^2 = \left[\sum_{1 \leq i < j \leq m}(p_i - p_j)^{-2}\right]^{-1}. \qquad (30)$$

It is obvious from the sum in the formula (30) is a symmetric function of the roots $p_1, \ldots, p_m$ of the polynomial (11). Therefore, according to the Theorem 1, it can be expressed as a rational function of the elementary symmetric functions (12) and, hence, as a rational function of the coefficients of the given polynomial $P_n(x)$. Then, to prove the statement of the Lemma for $\mu_k$, $k = 1, 2, \ldots$, let us note that the expression (17) which contain the function $y(\mu_0)$, is symmetric respect to permutations of the roots $p_1, \ldots, p_m$. As far as $\mu_0$ also is a symmetric function of these roots, the statement of the Lemma is valid for $\mu_1$. Then, using the mathematical induction one can proof the statement of the Lemma for any $\mu_k$, $k = 1, 2, \ldots$. Now, in order to proof the statement for the $M_k$, $0 < M_k < M$, $k = 0, 1, 2, \ldots$, it is enough to express $M_0$ as a rational function of the roots $p_1, \ldots, p_m$. One can do it applying the formula (3') to the polynomial $P_m(x)$:

$$D_2(P_m) = \sum_{1 \leq i < j \leq m}(p_j - p_i)^2 = s_0 s_2 - s_1^2 = m(\sigma_1^2 - 2\sigma_2) - \sigma_1^2 > M^2 = \max(p_j - p_i)^2.$$

Hence, one can choose

$$M_0^2 = D_2(P_m). \qquad (31)$$

According to the Theorem 1, it means that the Lemma is valid for $M_0$ and then, due to mathematical induction, for any $0 < M_k < M$, $k = 0, 1, 2, \ldots$. ∎

So, we have proven the next statement:

**Theorem 2**. The minimal (the maximal) distance between the (real) roots of the polynomial can be find out as a limit of an increasing (decreasing) convergent sequence of the rational functions which depend on the coefficients of the polynomial only. ∎

**Remark 1**. The function $Z(\varepsilon)$ defined by the formula (16) is presented in the subsequent consideration only as a combination $\frac{1}{\varepsilon}\left(Z(\varepsilon) + \frac{m}{2\varepsilon}\right) = \sum_{1 \leq i < j \leq m}\left[(p_i - p_j)^2 - \varepsilon^2\right]^{-1}$. So, one has to express by coefficients of given polynomial the last symmetric function only. This is much easier than to express the function $Z(\varepsilon)$ based on the formulas (13) and (14).

**Remark 2.** According to the Lemma, all terms of the convergent sequences $\{\mu_k | k=0,1,...\}$ and $\{M_k | k=0,1,...\}$ are symmetric functions of the roots $p_i$, $i=\overline{1,m}$, $m \geq 3$, of the polynomial. It seems interesting to investigate whether the limits of these sequences $\mu_\infty = \min |p_i - p_j|$ and $M_\infty = \max |p_i - p_j|$ are also symmetric functions of the roots $p_i$, $i = \overline{1,m}$.

**Remark 3.** According to the definition (16), for $m=1$ one has

$$Z(\varepsilon) = -\frac{1}{2\varepsilon} \quad \Rightarrow \quad \frac{1}{\varepsilon}\left(Z(\varepsilon) + \frac{1}{2\varepsilon}\right) = 0; \quad (m=1).$$

Correspondingly, for $m=2$ one gets

$$Z(\varepsilon) = -\frac{1}{\varepsilon} + \frac{\varepsilon}{(p_2 - p_1)^2 - \varepsilon^2} \quad \Rightarrow \quad \frac{1}{\varepsilon}\left(Z(\varepsilon) + \frac{1}{\varepsilon}\right) = \frac{1}{(p_2 - p_1)^2 - \varepsilon^2}. \quad (m=2).$$

According to the formulas (30), (31) and (3'), in the last case, obviously, we have

$$\mu = p_2 - p_1 = \mu_0 = M_0 = M.$$

Let us estimate the convergence rate of the sequences $\{\mu_k | k=0,1,...\}$ and $\{M_k | k=0,1,...\}$ ($m \geq 3$). Taking into account that for any $k = 0,1,...$ the inequality $\mu_k < \mu$ is valid, one gets for the first sequence:

$$\mu^2 - \mu_k^2 > \mu_{k+1}^2 - \mu_k^2 = \left(\sum_{1 \leq i < j \leq m} \frac{1}{(p_i - p_j)^2 - \mu_k^2}\right)^{-1} \geq \left(\sum_{l=1}^{m-1} \frac{m-l}{l^2\mu^2 - \mu_k^2}\right)^{-1}. \quad (32)$$

The sum in the right hand side of the formula (32) corresponds to a polynomial having equidistant roots with the minimal distance between them $\mu$. Let us denote this polynomial

$$W_m(x) = \prod_{l=0}^{m-1}(x - \mu l), \quad \mu > 0. \quad (33)$$

and call it "Wilkinson's generalized polynomial". The Wilkinson's polynomial [7] corresponds to the case $\mu = 1$.

Calculating the sum in (32) for the polynomial $W_m(x)$ one obtains

$$\sum_{l=1}^{m-1} \frac{m-l}{l^2\mu^2 - \mu_k^2} = \frac{m-1}{\mu^2 - \mu_k^2} + \sum_{l=2}^{m-1} \frac{m-l}{l^2\mu^2 - \mu_k^2} =$$

$$= \frac{1}{\mu_k^2}\left(\frac{m-1}{\mu^2/\mu_k^2 - 1} + \sum_{l=2}^{m-1} \frac{m-l}{l^2\mu^2/\mu_k^2 - 1}\right) < \frac{1}{\mu_k^2}\left(\frac{m-1}{\mu^2/\mu_k^2 - 1} + \sum_{l=2}^{m-1} \frac{m-l}{l^2 - 1}\right).$$

It is not difficult to calculate the last sum which, obviously, is positive ($m \geq 3$) and increasing when $m$ increases:

$$\sum_{l=2}^{m-1} \frac{m-l}{l^2 - 1} = \frac{3}{4}m - \frac{1}{4} + \frac{1}{2m} - \sum_{l=1}^{m} l^{-1} = \frac{3}{4}m - \frac{1}{4} + \frac{1}{2m} - \left[\frac{1}{m} + \gamma + \psi(m)\right] < \frac{3}{4}m.$$

The last inequality is valid as far as $\gamma + \psi(m) \geq 0$ for any $m \geq 1$. Here $\psi(z) = d/dz \ln \Gamma(z)$ denotes the digamma function, $\Gamma(z)$ being the Euler gamma function and $\gamma$ is the Euler-Mascheroni constant, $\gamma = -\psi(1) = 0.577...$ [6]. Hence, one obtains the estimations

$$0 < \left[\frac{m-1}{\mu^2/\mu_k^2 - 1} + \frac{3}{4}m\right]^{-1} < \frac{\mu_{k+1}^2}{\mu_k^2} - 1 < \frac{\mu^2}{\mu_k^2} - 1. \quad (34)$$

As far as $\mu_k^2 \underset{k\to\infty}{\to} \mu^2$, we can choose the iteration number $k \geq k^*$ according to the condition

$$0 < \frac{\mu^2 - \mu_k^2}{\mu_k^2} < \frac{4}{3m}, \quad k \geq k^*. \tag{35}$$

Then we can simplify the left boundary in the estimation (34) and get the new estimation:

$$\frac{4}{3m^2} < \frac{\mu_{k+1}^2}{\mu_k^2} - 1 < \frac{4}{3m}, \quad (m \geq 3, \ k \geq k^*). \tag{36}$$

Now let us estimate a convergence rate for the sequence $\{M_k \mid k = 0, 1, \ldots\}$ ($m \geq 3$). As far as for all $i, j = \overline{1, m}$ we have $(p_i - p_j)^2 \leq M^2 < M_k^2$, rather similarly to the previous case we get

$$0 < \left[ M_k^2 - M^2 \right]^{-1} < \left[ M_k^2 - M_{k+1}^2 \right]^{-1} = \sum_{1 \leq i < j \leq m} \frac{1}{M_k^2 - (p_i - p_j)^2} < \frac{m(m-1)}{2} \left[ M_k^2 - M^2 \right]^{-1}.$$

$$(m \geq 3)$$

and finally we obtain

$$0 < \frac{2}{m(m-1)} \left( 1 - \frac{M^2}{M_k^2} \right) < 1 - \frac{M_{k+1}^2}{M_k^2} < 1 - \frac{M^2}{M_k^2}. \quad (m \geq 3) \tag{37}$$

Taking the difference between the right and the left hand sides' boundaries equal to a precision $\delta$ of numerical computations one obtains

$$\frac{(m+1)(m-2)}{m(m-1)} \left( 1 - \frac{M^2}{M_k^2} \right) = \delta, \quad 1 - \frac{M^2}{M_k^2} = \frac{m(m-1)}{(m+1)(m-2)} \delta \tag{38}$$

and therefore

$$\frac{2}{(m+1)(m-2)} \delta < 1 - \frac{M_{k+1}^2}{M_k^2} < \left( 1 + \frac{2}{(m+1)(m-2)} \right) \delta. \tag{39}$$

The estimation (39) specifies the step $k$ where one should stop iterations in the recurrence formula (24).

**Theorem 3**. The segment of number axis where all roots of the polynomial (having real roots only) are located can be estimate by rational functions of the coefficients of the polynomial.

**Proof.** Let us estimate the maximum distance between the roots of the polynomial and the arithmetic mean of the roots $\overline{p} = \frac{1}{m} s_1$. Using the method developed above, let us define

$$\left( \sum_{1 \leq i \leq m} \frac{1}{(p_i - \overline{p})^2 - m_k^2} \right)^{-1} + m_k^2 = m_{k+1}^2, \quad k = 0, 1, 2, \ldots . \tag{40}$$

In order to obtain decreasing and lower bounded sequence:

$$(p_i - \overline{p})^2_{\max} < m_k^2 < \ldots < m_2^2 < m_1^2 < m_0^2, \tag{41}$$

one has to choose $(p_i - \overline{p})^2_{\max} < m_0^2$, that is provided by taking as $m_0$ any term from the sequence (23).

One can stop at any $k$ and consider the values $\overline{p} - \sqrt{m_k^2} \equiv a$ and $\overline{p} + \sqrt{m_k^2} \equiv b$. Obviously, all roots belong to the segment $[a, b]$. Generally, it is wider than $M = \max |p_i - p_j|$, therefore, $a, b \neq p_i$, $i = \overline{1, m}$.

Now we can find the minimal distance between the roots and ends of the segment $[a, b]$. Choosing

$$\alpha_0^2 = \left[ \sum_{1 \leq i \leq m} (p_i - a)^{-2} \right]^{-1} < \min(p_i - a)^2, \quad i = \overline{1, m}, \tag{42}$$

and defining $\alpha_k$ as

$$\alpha_{k+1}^2 \equiv \left( \sum_{1 \leq i \leq m} \frac{1}{(p_i - a)^2 - \alpha_k^2} \right)^{-1} + \alpha_k^2, \quad k = 0, 1, 2, \ldots, \tag{43}$$

one obtains increasing and over bounded sequence

$$\alpha_0^2 < \alpha_1^2 < \alpha_2^2 < \ldots < \alpha_k^2 < \min(p_i - a)^2, \quad i = \overline{1, m}, \tag{44}$$

having the finite limit

$$\alpha_\infty^2 = \min(p_i - a)^2, \quad i = \overline{1, m}. \tag{45}$$

Similarly, one constructs increasing and over bounded sequence converging to the minimal distance between the roots and upper edge of the segment $b$:

$$\beta_0^2 < \beta_1^2 < \beta_2^2 < \ldots < \beta_k^2 < \ldots < \beta_\infty^2 = \min(p_i - b)^2, \quad i = \overline{1, m}, \tag{46}$$

where

$$\beta_{k+1}^2 \equiv \left( \sum_{1 \leq i \leq m} \frac{1}{(p_i - b)^2 - \beta_k^2} \right)^{-1} + \beta_k^2, \quad k = 0, 1, 2, \ldots, \tag{47}$$

$$\beta_0^2 = \left[ \sum_{1 \leq i \leq m} (p_i - b)^{-2} \right]^{-1} < \min(p_i - b)^2, \quad i = \overline{1, m}. \tag{48}$$

It's clear that $p_i \in (a - \alpha_k, b - \beta_k)$, $i = \overline{1, m}$. Taking into account that all expressions in the formulas (40), (42), (43), (47), (48), as well as the arithmetic mean of the roots $\bar{p}$, are symmetric functions of the polynomial's roots, one gets the statement of the Theorem. ∎

For the Wilkinson's generalized polynomial (33) a lot of estimations can be performed explicitly. First of all, in this case the formula (30) gives us

$$\mu_0^2 = \left[ \sum_{1 \leq i < j \leq m} \mu^{-2}(i-j)^{-2} \right]^{-1} = \mu^2 \left( \sum_{l=1}^{m-1} \frac{m-l}{l^2} \right)^{-1} = \mu^2 w_0^2(m),$$

where we have denoted

$$w_0^2(m) = \left[ m \frac{\pi^2}{6} - m\psi'(m) - \gamma - \psi(m) \right]^{-1} > 0. \quad (m \geq 3) \tag{49}$$

**Statement 2**. For polynomial (33) and for any $k = 0, 1, \ldots$ one has

$$\mu_k^2 = \mu^2 w_k^2(m), \quad k = 0, 1, 2\ldots, \quad (m \geq 3) \tag{50}$$

and the coefficients $w_k^2(m)$ satisfy the recurrence relation

$$w_{k+1}^2(m) = w_k^2(m) + \left( \sum_{l=1}^{m-1} \frac{m-l}{l^2 - w_k^2(m)} \right)^{-1}, \quad k = 0, 1, 2\ldots . \tag{51}$$

**Proof.** By induction. ∎

**Statement 3**. One has

$$w_0^2 < w_1^2 < \ldots < w_k^2 < 1, \tag{52'}$$

$$\lim_{k \to \infty} w_k(m) = 1, \tag{52''}$$

and, as far $w_k(m) = 1 - \varepsilon$ $(\varepsilon > 0)$, the convergence rate of the recurrence relation (51) can be estimated as

$$\frac{4\varepsilon}{7m - 4} < w_{k+1}^2(m) - w_k^2(m) < \frac{1}{m-1}. \tag{53}$$

**Proof.** Relations (52′), (52″) can be easily proven similarly to the general formulae (19) and (21). The estimation (53) can be proven using the condition $w_k(m) = 1 - \varepsilon$, $0 < \varepsilon < 1$ and presenting the relation (51) as:

$$w_{k+1}^2(m) - w_k^2(m) = \left(\sum_{l=1}^{m-1} \frac{m-l}{l^2 - (1-\varepsilon)^2}\right)^{-1} = \left(\frac{m-1}{\varepsilon(2-\varepsilon)} + \tilde{w}(m)\right)^{-1},$$

where

$$\tilde{w}(m) \equiv \sum_{l=2}^{m-1} \frac{m-l}{l^2 - (1-\varepsilon)^2}.$$

Due to inequality $0 < \varepsilon < 1$ one has

$$\sum_{l=2}^{m-1} \frac{m-l}{l^2} = w_0^{-2}(m) - (m-1) < \tilde{w}(m) < \sum_{l=2}^{m-1} \frac{m-l}{l^2 - 1}$$

and therefore

$$\left(\frac{m-1}{\varepsilon(2-\varepsilon)} + \sum_{l=2}^{m-1} \frac{m-l}{l^2 - 1}\right)^{-1} < \left(\frac{m-1}{\varepsilon(2-\varepsilon)} + \tilde{w}(m)\right)^{-1} = w_{k+1}^2(m) - w_k^2(m) < \left(\frac{m-1}{\varepsilon(2-\varepsilon)} + \sum_{l=2}^{m-1} \frac{m-l}{l^2}\right)^{-1}.$$

Taking here into account that $0 < \varepsilon(2-\varepsilon) < 1$, $0 < \varepsilon(2-\varepsilon) < 2\varepsilon$, and, besides, that the both finite sums are positive ($m \geq 3$; see the formula (49)):

$$0 < \sum_{l=2}^{m-1} \frac{m-l}{l^2 - 1} = \frac{3}{4}m - \frac{1}{4} - \frac{1}{2m} - (\gamma + \psi(m)) < \frac{3}{4}m, \quad 0 < \sum_{l=2}^{m-1} \frac{m-l}{l^2} = w_0^{-2} - (m-1),$$

one obtains the estimations

$$\frac{4\varepsilon(2-\varepsilon)}{4(m-1) + 3m\varepsilon(2-\varepsilon)} < w_{k+1}^2(m) - w_k^2(m) < w_0^2(m) < \frac{1}{m-1},$$

which can be simplified, using easily checkable inequalities:

$$\frac{4\varepsilon(2-\varepsilon)}{4(m-1) + 3m\varepsilon(2-\varepsilon)} > \frac{4\varepsilon(2-\varepsilon)}{4(m-1) + 3m} = \frac{4\varepsilon(1+1-\varepsilon)}{7m-4} > \frac{4\varepsilon}{7m-4}.$$

So, we get the estimations (53). ∎

**Corollary 4.** Obviously,

$$w_{k+1}^2(m) - w_k^2(m) = (w_{k+1}(m) + w_k(m))(w_{k+1}(m) - w_k(m)) < 2(w_{k+1}(m) - w_k(m))$$

and then, according to the estimation (53) one obtains

$$w_{k+1}(m) - w_k(m) > \frac{1}{2}(w_{k+1}^2(m) - w_k^2(m)) > \frac{2\varepsilon}{7m-4}.$$

Hence,

$$w_{k+1}(m) > w_k(m) + \frac{2\varepsilon}{7m-4}. \tag{54}$$

The last inequality is valid for any $k = 0, 1, \ldots$ and any $0 < \varepsilon < 1$. Obviously, $\varepsilon$ depends on $k$ and $m$:

$$\varepsilon = 1 - w_k(m) \equiv \varepsilon_k(m).$$

Hence, the estimation (54) gives

$$\varepsilon_k(m) - \varepsilon_{k+1}(m) > \frac{2\varepsilon}{7m-4}.$$

Supposing here $\varepsilon = \varepsilon_k(m)$ we obtain the next recurrence formula:

$$\varepsilon_{k+1}(m) < \varepsilon_k(m) \frac{7m-6}{7m-4}, \qquad k = 0,1,\ldots.$$

Therefore, the sequence $\{\varepsilon_k(m) \mid k = 0,1,\ldots\}$ is majored (dominated) by an exponential sequence - converging geometric progression

$$\{b_k \mid b_k = \varepsilon_0 q^k, k = 0,1,\ldots\},$$

with exponentially denominator

$$q(m) = \frac{7m-6}{7m-4} = 1 - \frac{2}{7m-4} < 1.$$

Hence, one has

$$\varepsilon_{k+1}(m) < \varepsilon_0(m) q^k(m) \approx \varepsilon_0(m)\left(1 - \frac{2k}{7m-4}\right), \tag{55}$$

and taking into account that, according to (52'),

$$0 < w_0 < 1,$$

$$\varepsilon_0(m) = 1 - w_0 = 1 - \left[m\frac{\pi^2}{6} - m\psi'(m) - \gamma - \psi(m)\right]^{-1/2}, \tag{56}$$

$$0 < \varepsilon_0(m) < 1, \quad (m \geq 3)$$

one concludes that $\varepsilon_k$ becomes about the precision of numerical calculation $\delta$ for

$$k \approx \frac{7m-4}{2}(1 - \delta/\varepsilon_0). \tag{57}$$

## Conclusion

The coefficients of the minimal polynomial of the Hermitian matrix are rationally expressed by the coefficients of the characteristic polynomial. There are find out (with any precision) analytic estimations for minimal and maximal distances between roots of the minimal polynomial as well as the range of localization of all its roots. Factorizing the given characteristic polynomial in such manner that all roots with the certain multiplicity belong to the certain factor, we apply the method developed to each factor separately. This gives us analytic conditions for classification polynomials with respect to relative location of roots with certain multiplicities.